\documentclass[hidelinks]{article}
\newcommand{\documentdate}{30 I 2026}

\usepackage{a4wide,latexsym,amsmath}
\usepackage{varioref}
\usepackage{xcolor}
\usepackage{amsfonts}
\usepackage{hyperref}
\usepackage{todonotes}
\usepackage{bbm}
\title{Fast Stochastic AdaGrad\\
for Nonconvex Bound-Constrained Optimization \\ using Second-Order Information}

\author{Stefania Bellavia\footnotemark[1],
        Serge Gratton\footnotemark[2],
        Benedetta Morini\footnotemark[3],
        Philippe L. Toint\footnotemark[4]}

\newcommand{\beqn}[1]{\begin{equation}\label{#1}}
\newcommand{\eeqn}{\end{equation}}
\newcommand{\req}[1]{(\ref{#1})}
\newcommand{\ms}{\;\;\;\;}
\newcommand{\tim}[1]{\;\; \mbox{#1} \;\;}

\newcounter{algo}[section]
\renewcommand{\thealgo}{\thesection.\arabic{algo}}
\newcommand{\llem}[2]{\vspace{\baselineskip} 
\noindent\framebox[\textwidth]{\parbox{0.95\textwidth}{
\begin{lemma} \label{#1} \rm #2 \end{lemma} } } \vspace{\baselineskip} }
\newcommand{\algo}[3]{\refstepcounter{algo}
\begin{center}\begin{figure}[htbp]
\framebox[1.05\textwidth]{
\parbox{\textwidth} {\vspace{\topsep}
{\bf Algorithm \thealgo : #2}\label{#1}\\
\vspace*{-\topsep} \mbox{ }\\
{#3} \vspace{\topsep} }}
\end{figure}\end{center}}
\newcommand{\bpr}{{\bf Proof.} \hspace{1.5mm}}
\newcommand{\epr}{\hfill $\Box$ \vspace*{1em}}
\newcommand{\lthm}[2]{\vspace{\baselineskip} 
\noindent\framebox[\textwidth]{\parbox{0.95\textwidth}{
\begin{theorem} \label{#1} \rm #2 \end{theorem} } } \vspace{\baselineskip} }
\newcommand{\lcor}[2]{\vspace{\baselineskip} 
\noindent\framebox[\textwidth]{\parbox{0.95\textwidth}{
\begin{corollary} \label{#1} \rm #2 \end{corollary} } } \vspace{\baselineskip}
}
\newcommand{\ii}[1]{\{ 1, \ldots, #1 \}}
\newcommand{\iiz}[1]{\{ 0, \ldots, #1 \}}

\newcommand{\calB}{{\cal B}}

\newcommand{\calF}{{\cal F}} 
\newcommand{\calO}{{\cal O}} 

\newcommand{\calT}{{\cal T}}
\renewcommand{\Re}{\hbox{I\hskip -2pt R}}

\newcommand{\bigfrac}[2]{\frac{\displaystyle #1}{\displaystyle #2}}

\newcommand{\bigmax}{\displaystyle \max}
\newcommand{\sfrac}[2]{{\scriptstyle \frac{#1}{#2}}}
\newcommand{\half}{\sfrac{1}{2}}
\newcommand{\eqdef}{\stackrel{\rm def}{=}}
\newcommand{\bigsum}{\displaystyle \sum}

\newcommand{\kap}[1]{\kappa_{\mbox{\tiny #1}}}
\newcommand{\flow}{f_{\rm low}}
\newcommand{\al}[1]{{\footnotesize{\sf #1}}}
\newcommand{\tal}[1]{{\normalsize {\sf #1}}}

\newtheorem{theorem}{Theorem}[section]
\newtheorem{lemma}[theorem]{Lemma}
\newtheorem{corollary}[theorem]{Corollary}

\newcommand{\proof}[1]{
\begin{list}{}{
\setlength{\topsep}{0.0pt}
\setlength{\partopsep}{0.0pt}
\setlength{\leftmargin}{0.025\textwidth}
\setlength{\rightmargin}{0.5\leftmargin}
\setlength{\labelwidth}{0.5\leftmargin}
\setlength{\labelsep}{0.25\leftmargin}}
\item \bpr #1 \epr \noindent
\end{list}}
\DeclareMathOperator*{\average}{average}
\DeclareMathOperator*{\argmax}{argmax}
\newcommand{\E}[1]{\mathbb{E}\!\left[#1 \right]}

\newcommand{\Econd}[2]{\mathbb{E}_{#1}\!\left[#2 \right]}

\newcommand{\Prob}[1]{\mathbb{P}\!\left[#1 \right]}

\newcommand{\RMSE}{\mbox{\footnotesize RMSE}}
\newcommand{\new}[1]{{\color{blue}#1}}
\newcommand{\comment}[1]{}
\definecolor{privategreen}{rgb}{0,0.6,0}
\newcommand{\private}[1]{}

\newcommand{\algname}{\al{ADAGB2}}

\date{\documentdate}

\begin{document}

\maketitle

\renewcommand{\thefootnote}{\fnsymbol{footnote}}
\footnotetext[1]{Dipartimento  di Ingegneria Industriale, 
  Universit\`a degli Studi di Firenze,
  Firenze,  Italia. Member of the INdAM Research Group GNCS.
  Email: stefania.bellavia@unifi.it.}
\footnotetext[2]{Universit\'e de Toulouse, INP, IRIT, Toulouse,
  France. Work partially supported by the Artificial and Natural
  Intelligence Toulouse Institute (ANITI), French ``Investing for the
  Future - PIA3'' program under the Grant agreement ANR-19-PI3A-0004.
  Email: serge.gratton@enseeiht.fr.}
\footnotetext[3]{Dipartimento  di Ingegneria Industriale, 
  Universit\`a degli Studi di Firenze,
  Firenze,  Italia. Member of the INdAM Research Group GNCS. 
  Email: benedetta.morini@unifi.it.}
\footnotetext[4]{Namur Center for Complex Systems (naXys),
  University of Namur, Namur, Belgium.Work partially supported 
  by the  Artificial and Natural Intelligence Toulouse Institute (ANITI).
  Email: philippe.toint@unamur.be.}
\renewcommand{\thefootnote}{\arabic{footnote}}

\begin{abstract}
\algname, a generalization of the \tal{AdaGrad} algorithm for 
stochastic optimization is introduced, which handles
 bound-constrained problems and uses second-order 
information when available. It is shown that, given $\delta\in(0,1)$ and $\epsilon\in(0,1]$, the \algname\ algorithm needs at most 
$\calO(\epsilon^{-2})$ iterations to ensure an $\epsilon$-approximate 
first-order critical point of the bound-constrained problem with probability 
at least $1-\delta$ under a new ``directional'' assumption on the  gradient oracle and its conditional root-mean-square error. Should this latter condition fail, it is also 
shown that the optimality level of iterates is bounded above by this average.
The relation between the approximate and true classical projected-gradient-based optimality measures for bound constrained problems is also investigated, and it is shown that merely assuming unbiased gradient oracles may be insufficient to ensure convergence in $\calO(\epsilon^{-2})$ iterations.
\end{abstract}


\section{Introduction}

First-order optimization algorithms have been widely used in the
contexts on online learning and deep neural network training and their
convergence properties on nonconvex problems have been investigated by several authors
(see \cite{BottCurtNoce18} for a survey).  Among them, \tal{AdaGrad} 
\cite{DuchHazaSing11}, despite not being always as
efficient as others in practice on nonconvex problems \cite{ZhanXu12}, enjoys a special
position from the theoretical point of view because of its solid and extensive
convergence analysis.

When the objective function's gradient is contaminated by noise (for instance caused by sampling)
a probabilistic point of view on the algorithm's convergence theory is desirable.
This has been investigated for \tal{AdaGrad} applied to nonconvex functions by a number 
of authors, as shown in Table~\ref{table:adag-theory}, along with the relevant important 
assumptions and results. We discuss this rich body of theory (and the content of this table) 
in more detail in Section~\ref{section:survey}.

\begin{table}[thbp]
\begin{center}
\footnotesize
\begin{tabular}{lcccccccc}
Paper & & Smooth & Gradients & Noise & 2nd-order & Bound & Conv.& Conv. \\
      & &        & \hspace*{-3mm}bound~/~bias &  type &   usage  & constr. & type &rate \\
\hline
Li-Orabona \hspace*{-3mm} & \cite{LiOrab19} & $L^*$ & no~/~no & Sub-Gauss. & no & no & $\E{.}$/w.h.p. & $\calO_\ell(\epsilon^{-2})$\\
&  & &  & + $\sigma\searrow 0$ &  &  &   &\\
Ward et al. & \cite{ward2020adagrad} & $L$ &  yes~/~no  & Bounded & no & no & w.h.p. & $\calO_\ell(\epsilon^{-4})$\\
D\'{e}fossez et al. \hspace*{-3mm}& \cite{DefoBottBachUsun22} & $L$ & yes~/~no & Unrestricted & no & no & $\E{.}$ & $\calO_\ell(\epsilon^{-4})$\\
Gratton et al. \hspace*{-3mm}& \cite{GratJeraToin22b} & $L$ & yes~/~no & Unrestricted &  no & no & $\E{.}$  & $\calO_\ell(\epsilon^{-4})$\\
Kavis et al. \hspace*{-3mm}& \cite{KaviLevyCevh22} &$L$ & yes~/~no & Sub-Gauss. &  no & no & w.h.p  & $\calO_\ell(\epsilon^{-2})$\\
&  & &  & + $\sigma\searrow 0$ &  &  &   &\\
Faw et al. \hspace*{-3mm}& \cite{faw2022power}& $L$ & no~/~no & Affine & no &  no & w.h.p.  & $\calO_\ell(\epsilon^{-4})$\\
Wang et al. \hspace*{-3mm}& \cite{wang2023convergence} & $(L_0,L_1^*)$ & no~/~no & Affine & no & no & $\E{.}$/w.h.p.  & $\calO_\ell(\epsilon^{-2})$\\
&  & &  & + $\sigma\searrow 0$ &  &  &   &\\
Liu et al. \hspace*{-3mm}& \cite{Liuetal23} & $L$ & no~/~no & Sub-Gauss. & no & no & w.h.p. & $\calO(\epsilon^{-2})$\\
 & & &  & + $\sigma\searrow 0$ &  &  &   &\\
Attia-Koren \hspace*{-3mm}& \cite{AttiKore23} & $L$ & no~/~no & Affine-* & no & no & w.h.p. & $\calO_\ell(\epsilon^{-2})$\\
&  & &  & + $\sigma\searrow 0$ &  &  &   &\\ 
Faw et al. \hspace*{-3mm}& \cite{faw2023beyond} & $(L_0,L_1^*)$  & no~/~no & Affine & no & no & w.h.p. & $\calO_\ell(\epsilon^{-4})$\\
Hong-Lin \hspace*{-3mm}& \cite{HongLin24} & $L$ & no~/~no & Affine-* &  no  &  no & w.h.p. & $\calO_\ell(\epsilon^{-2})$\\
&  & &  & + $\sigma\searrow 0$ &  &  &   &\\ 
Hong-Lin \hspace*{-3mm}& \cite{HongLin24} & $(L_0,L_1^*)$ & no~/~no & Affine-* &  no  &  no & w.h.p. & $\calO_\ell(\epsilon^{-4})$\\
Jiang et al. \hspace*{-3mm}& \cite{JianMalaMokt24} & $L$ & no~/~no & Bounded &  no & no & $\E{.}$ & $\calO_\ell(\epsilon^{-2})$\\
Alacaoglu et al.\hspace*{-2mm}.& \cite{AlacMaliChev21,AlacLyu23} & L & yes/yes & Bounded & no & yes & $\E{.}$ & $\calO_M(\epsilon^{-4})$\\
This paper \hspace*{-3mm}& & $L$ &   no~/~no  & new1      &  yes & yes & $\E{.}$/w.h.p & $\calO(\epsilon^{-2})$\\
This paper \hspace*{-3mm}& & $L$ &   no~/~yes & new1+new2 &  yes & yes & $\E{.}$/w.h.p & $\calO(\epsilon^{-2})$\\
\end{tabular}
\caption{\label{table:adag-theory}\footnotesize Convergence theories for the stochastic \al{AdaGrad} algorithm and their characteristics in the nonconvex setting. 
 (In the third column, $L^*$ means that the Lipschitz constant must be known, the last column reports a bound (in order) for the considered algorithm to achieve $\epsilon$-criticality, the subscript $\ell$ indicating the presence of an additional logarithmic factor in $\epsilon$, the subscript $M$ indicates that the convergence is studied in terms of  the norm of the gradient of the Moreau envelope,
 see Section \ref{section:survey} for more detail on the other columns).}
\end{center}
\end{table}

These proposals are however limited, from the theoretician's
point of view, in three respects.  
\begin{itemize}
\item 
The first is that the probabilistic analysis is restricted to the case 
where second-order information is ignored when available. The algorithm 
remains strictly first-order, with a step always aligned (possibly component-wise) 
with the approximate steepest-descent direction (no trust region is used).
By contrast, other first-order methods have been "augmented" to use step-sizes along the 
first-order direction taking second-order information into account (see 
\cite{BelMorYou2024, BerJahRicTak2022,CurShi2022,FarAssMur2023, Tanetal16,WanWuMat2021} and the 
references therein). 
\item
The second is the capability of solving constrained problems that arise 
 when \textit{a priori} information on the
problem at hand is available, often in the form of
constraints, of which bounds on the variables are the most common. 
Such problems arise, for instance in fracture mechanics 
\cite{kopanivcakova2020recursive, krause2009nonsmooth}, inverse problems and 
PDE-constrained optimization under uncertainty \cite{ciaramella2024multigrid}. 
Moreover, bounds are frequently applied to machine learning models to avoid overfitting \cite{leimkuhler2021better} or to reflect real-world limits
and maintain physical interpretability~\cite{lu2021physics}.
Admittedly, these can be taken care of using
penalty terms in the loss/objective function, as is often done
for Physically Informed Neural Networks (PINNs)
(see \cite{CaiMaoWangYinKarn21,WangYuPerd22} for instance), but this
introduces new hyper-parameters needing calibration and does not
ensure that constraints are strictly satisfied. Moreover, a
strong penalization of the constraints degrades
the problem's conditioning, possibly causing slow
convergence, especially if first-order methods are used.

Focusing on  \tal{Adagrad} for constrained problems, 
we are aware of only \cite{AlacLyu23,AlacMaliChev21} that provide the stochastic analysis of
 a projected version of  \tal{Adagrad-Norm}. 
 More generally, stochastic gradient-based methods for constrained problems have been proposed, e.g., in \cite{AlacLyu23,CurtKungRobiWang23,DaviDrus19,ghadimi2016composite,Wangpierzhoucurt26} and references there-in,
but these techniques differ considerably from the \tal{AdaGrad}-like gradient algorithms considered here. 

\item 
The third shortcoming is that, to the best of the authors' knowledge, bias in the 
gradient oracle has been considered only for \tal{AdaGrad-Norm} in \cite{AlacLyu23,AlacMaliChev21}, in contrast with other 
stochastic optimization techniques (see \cite{BeraCaoSche21, BlanCartMeniSche19} for 
instance).
\end{itemize}

\noindent
The challenge addressed in this paper is thus threefold. We first consider a 
probabilistic theory of how second-order steps can be accommodated in a stochastic 
\tal{AdaGrad}-like\footnote{That is an algorithm that reduces to \tal{AdaGrad} when 
second-order information is not used.} optimization method. We also show how this 
can be done in the context of bound-constrained problems, and, finally, propose a 
theory which does not assume that gradient oracles are unbiased.

More specifically, we consider the problem
\beqn{problem}
\min_{x \in \calF} f(x)
\tim{ where }
\calF = \{x\in\Re^n \mid l_i \leq x_i \leq u_i \}
\eeqn
where $f$ is a smooth (possibly nonconvex) function from an open set containing the 
feasible region $\calF\subseteq \Re^n$ into $\Re$, and where $l$ and $u$ are vectors
specifying the lower and upper bounds on the variables, respectively (infinite entries
in $l$ and/or $u$ are allowed). We also assume that $\nabla_x^1f(x)$ cannot be computed 
exactly but is approximated by a random oracle.  As a consequence, the 
algorithm we are about to describe generates a random process, where, for a given 
iterate $x_k$, the oracle computes the gradient oracle $g(x_k,\xi)$ where $\xi$ is a 
random variable (whose distribution may depend on $x_k$), with probability 
space $(\Omega, \calT_\Omega, \mathbb{P})$. For brevity, we will denote $g_k\eqdef 
g(x_k,\xi)$ and we also define $G_k = \nabla_x^1f(x_k)$ and $H_k = \nabla_x^2f(x_k)$.
 Conditioned to knowing $g_0,\ldots,g_{k-1}$, the expectation will be denoted by the 
symbol $\Econd{k}{\cdot}$, the probability of an event will be denoted as $P_k(\cdots)$ %\mathbb P_k({\rm event})$,
and $\mathbbm{1}(\cdot)$ will denote the indicator function of the event.
The symbol $\|\cdot\|$ denotes the 2-norm.

\subsection{Related works}\label{section:survey}

The convergence of  \tal{AdaGrad} in the nonconvex setting has been studied by a 
number of authors assuming the use of an unbiased stochastic oracle of the true 
gradient with a variety of assumptions on its noise.
In \cite{LiOrab19}, a  version of AdaGrad with delayed step-sizes was introduced. The  almost sure asymptotic convergence of the gradients to zero 
was proved for this variant assuming  that the noise has bounded support, i.e.,
\[
\Econd{k}{\|g_k-G_k\|} \leq \kap{bounded} \quad \quad  \kap{bounded} >0,
\]
at any iteration $k$. Assuming instead that, at each iteration $k$,
the sub-Gaussian noise condition
\[
\Econd{k}{e^{\|g_k-G_k\|^2/\sigma^2}} \leq e,
\]
holds, Li and Orabona also showed that \tal{AdaGrad} is noise adaptive, in 
the sense that, given $\epsilon>0$, the iteration complexity  interpolates,  with high 
probability (denoted by w.h.p. in Table~\ref{table:adag-theory}), between 
$\calO_\ell(\epsilon^{-4})$ to  $\calO_\ell(\epsilon^{-2})$ 
depending on the magnitude of $\sigma$, where $\calO_\ell(\cdot)$ means that 
the order is up to a logarithmic term.
Unfortunately, this analysis requires \textit{a priori} knowledge of the Lipschitz 
constant of the gradient in setting the step-size. A similar result was proved in 
\cite{KaviLevyCevh22} but without assuming this knowledge. 
 Ward, Wu and Bottou in \cite{ward2020adagrad} analyzed the convergence of 
 \tal{AdaGrad-Norm} assuming uniformly bounded gradients and  
 bounded variance of the gradient oracle at each iteration $k$, that is
 \[
\Econd{k}{\|g_k-G_k\|^2} \leq \sigma^2.
\]
They showed that the average of the squared norm of the gradients produced
by \tal{AdaGrad-Norm} converges with high probability at the rate 
$\calO(1/\sqrt{k})$, which implies that the algorithm needs at most 
$\calO_\ell(\epsilon^{-4})$ iterations to achieve $\E{\|G_k\|}\leq\epsilon$. 
Under weaker assumptions, i.e., without  assuming boundedness of the 
variance of the gradient's oracles, these  results have been extended  
(now in expectation form, noted by $\E{\cdot}$ in Table~\ref{table:adag-theory})
to the component-wise version of \tal{AdaGrad} in \cite{DefoBottBachUsun22}
and to an extended class of methods comprising \tal{AdaGrad} in \cite{GratJeraToin22b}.

A different stream of research analyzed \tal{AdaGrad} as an optimally-tuned 
adaptive stochastic gradient descent without assuming bounded gradients 
but under the affine bound on the variance at each iteration given by
\beqn{affine}
\Econd{k}{\|g_k-G_k\|^2} \leq \kap{1,affine} + \kap{2,affine} \|G_k\|^2,
\eeqn
with $\kap{1,affine}\ge 0$ and $\kap{2,affine} \ge 0$.  
Under these assumptions, Faw et al.~\cite{faw2022power} 
have shown that \tal{AdaGrad-Norm} exhibits an iteration  complexity of 
the order of $\calO_\ell(\epsilon^{-4})$, with high probability.

The noise adaptivity of \tal{AdaGrad} under the stronger affine condition
\beqn{affine*}
\|g_k-G_k\|^2 \stackrel{as}{\leq} \kap{1,affine} + \kap{2,affine} \|G_k\|^2
\eeqn
(denoted Affine-* in Table \ref{table:adag-theory}) has been proved in 
\cite{AttiKore23}. The same complexity bound for \tal{AdaGrad-Norm} 
is obtained in \cite{wang2023convergence} under the affine condition \req{affine}.

Very recently, Jiang, Maladkar and Mokhtari  in \cite{JianMalaMokt24} performed 
the convergence analysis in 1-norm and proved $\calO_\ell(\epsilon^{-2})$ 
iteration complexity, assuming coordinate-wise bounded variance and 
coordinate-wise Lipschitz continuity of the gradient.

Finally Hong and Lin revisited the convergence of \tal{AdaGrad} in the recent paper 
\cite{HongLin24}, assuming a relaxed version of the condition \req{affine*} given by
\beqn{affine*ext}
\|g_k-G_k\|^2 
\,\stackrel{as}{\leq}\, 
\kap{1,affine} + \kap{2,affine} \|G_k\|^2+ \kap{3,affine} (f(x_k)- \flow),
\eeqn
where $\flow$ is such that $f(x)\ge \flow\; \forall x \in \Re^n$. They prove 
that the iteration complexity interpolates between $\calO_\ell(\epsilon^{-4})$ 
to  $\calO_\ell(\epsilon^{-2})$ depending on the magnitude of $\kap{1,affine}$
and $\kap{3,affine}$. 

Projected versions of \tal{AdaGrad-Norm} have been
proposed in \cite{AlacMaliChev21,AlacLyu23}. In these latter papers biased stochastic oracles of the gradient are allowed. 
The paper \cite{AlacMaliChev21} elaborates on \cite{DaviDrus19} and proves the $\calO(\epsilon^{-4})$
iteration complexity employing as optimality measure 
 the norm of the gradient of Moreau envelope.  The  analysis works with mini-batch size of $1$ and assuming bounded stochastic gradients. This analysis has been further extended  under general non-i.i.d. data sampling assumption
 in \cite{AlacLyu23}.

All the previous mentioned results have been obtained assuming the Lipschitz 
continuity of the gradient. The papers  
\cite{faw2023beyond, HongLin24, wang2023convergence} provide results for the 
class of $(L_0,L_1)$-smooth functions\footnote{A function is 
$(L_0,L_1)$-smooth if there exist a constant $L_0 \geq0$ 
and $L_1 \geq 0$ and $\delta>0$ such that for all 
$x,y \in \Re^n$ with $\|x-y\|\leq \delta$, one has that 
$\|\nabla_x^1f(x)-\nabla_x^1f(y)\| \leq (L_0+L_1 \|\nabla_x^1f(y)\|)\|x-y\|$. } 
and assume that the constant $L_1$ is known. In \cite{faw2023beyond} the results 
in \cite{faw2022power} have been extended  to this latter  class of  functions and 
it has been proved that  the method  has a complexity of the order of 
$\calO_\ell(\epsilon^{-4})$, assuming \req{affine} with  $\kap{2,affine} <1$. 
Similar results are given in \cite{HongLin24} for \tal{AdaGrad} with momentum 
under the Assumption \req{affine*ext}. Stronger results are obtained in 
\cite{wang2023convergence} for \tal{AdaGrad-Norm} assuming 
\req{affine} where it is proved that the iteration complexity  
interpolates between $\calO_\ell(\epsilon^{-4})$ to  $\calO_\ell(\epsilon^{-2})$ 
depending on the magnitude of $\kap{1,affine}$.

To conclude this brief survey, we also note that Liu et al. \cite{Liuetal23}
have obtained an $\calO(\epsilon^{-2})$ convergence result for the 
objective-function gap produced by the stochastic \tal{AdaGrad-Norm} 
when applied to $\gamma$-quasar convex $L$-smooth functions with a sub-Weibull 
assumption on the gradient error. This condition subsumes the sub-Gaussian 
case, but the result does not apply to general $L$-smooth nonconvex functions.

\subsection{Summary of contributions}

In view of all contributions discussed above, we summarize our contributions as follows.
\begin{enumerate}
\item We propose an extension of \tal{AdaGrad} which is capable of using (possibly very approximate) second-order information whenever available.
\item This algorithm is also suitable for problems involving bounds 
on the variables. 
\item We analyze the convergence and probabilistic complexity of 
this algorithm, showing that it solves the approximate minimization 
problem with optimal complexity under a new directional condition 
on the gradient error, but without assuming bounded or unbiased 
gradient oracles nor knowledge of the problem's Lipschitz constant.
\item Taking a more general perspective, we discuss the 
relation between the classical optimality 
measure for bound-constrained problems induced by the projected gradient in 
the stochastic case, and have shown that, in general, the unbiased 
nature of the gradient oracle is not sufficient to ensure convergence 
on the true problem with the optimal rate, even if such a convergence 
occurs for the approximate one.

\item
We furthermore show that our optimal complexity result 
for the approximated problem extends to the true
problem (i.e.\ using exact gradients) if another condition on 
the gradient noise holds or if the gradient oracle is unbiased 
and the problem is unconstrained. We also describe the behaviour 
of the algorithm if none of these conditions hold.
\item
As far as the authors are aware, 
this is the   first convergence-rate result for the 
stochastic \tal{AdaGrad} applied to bound constrained problems and allowing to use second-order information. 
\end{enumerate}
\noindent
Section~\ref{section:algo} presents the new algorithm and discusses its features, while convergence analysis for the approximate problem is described in Section~\ref{section:convergence}. 
Section~\ref{section:spg} discusses the relation between approximate and true optimality measures in the bound-constrained case and its application to the new algorithm. A brief conclusion and some perspectives are finally outlined in Section~\ref{section:conclusion}.

\section{The algorithm}\label{section:algo}

Our proposed algorithm, called \algname, is presented~\vpageref{the-algo}.
Its inputs consist in a starting point $x_{\rm ini}$ and values $l$ and $u$ for the lower and upper bounds on the variables.

\algo{the-algo}{\algname
    ($x_{\rm ini},l,u$)}{
\begin{description}
\item[Step 0: Initialization:] The constants
$
\varsigma, \tau \in (0,1]$ and $\kappa_s \geq 1
$
are given.\\
Set  $x_0 = P_\calF(x_{\rm ini})$, $k=0$ and $w_{-1,i} =\sqrt{\varsigma}$ for $i\in\ii{n}$. 
\item[Step 1: First-order step:]
  Compute $g_k= g(x_k,\xi)$ a random approximation of $G_k$, set
  \beqn{dk-def}
   d_k \eqdef P_\calF(x_k - g_k ) - x_k,
  \eeqn
  \beqn{wDk-def}
  w_{k,i} = \sqrt{w_{k-1,i}^2+ d_{k,i}^2}
  \tim{and}
  \Delta_{k,i} =\frac{|d_{k,i}|}{w_{k,i}}
  \tim{ for } i\in \ii{n}.
  \eeqn
  \beqn{Bk-def}
  \calB_k
  = \left\{ x \in \Re^n \mid | x_i - x_{k,i}| \leq \Delta_{k,i}
  \tim{for} i \in\ii{n} \right\}
  \eeqn
  and
  \beqn{sL-def}
  s_k^L = P_{\calF\cap\calB_k}(x_k-g_k)-x_k.
  \eeqn
\item[Step 2: Second-order step: ] 
  Choose $B_k$ a symmetric approximation of $H_k$ and compute
  \beqn{gamma-def}
  s_k^Q = \gamma_ks_k^L
  \tim{where}
  \gamma_k = \left\{\begin{array}{ll}
      \min\left[ 1,\bigfrac{-g_k^Ts_k^L}{(s_k^L)^T B_k s_k^L}\right]
      & \tim{if } (s_k^L)^TB_ks_k^L > 0 \\
      1 & \tim{otherwise.}
    \end{array}\right.
  \eeqn
  Then select $s_k$ such that, for all $i\in \ii{n},$
  \beqn{inside-and-decr}
  x_k+s_k \in \calF,
  \ms
  |s_{k,i}|\leq \kappa_s\Delta_{k,i}
  \tim{ and }
  g_k^Ts_k+\half s_k^TB s_k
  \leq \tau\left(g_k^Ts_k^Q+\frac 1 2 (s_k^Q)^TB
  s_k^Q\right).
  \eeqn
\item[Step 3: Loop:] Set $x_{k+1} = x_k + s_k$, 
  increment $k$ by one and go to Step~1.
\end{description} 
}
\newpage

The statement of the algorithm suggests the following comments.
\begin{enumerate}
\item 
We first note that the mechanism of the algorithm ensures that all iterates remain feasible, that is $x_k\in \calF$ for all $k\geq 0$.
\item The projections $P_\calF$ and $P_{\calF\cap\calB_k}$ occurring in \req{dk-def} and \req{sL-def} are
extremely cheap to compute component-wise, as, for any vector $y\in
  \Re^n$ and $i\in\ii{n}$,
  \beqn{cheap-projs1}
  [P_\calF(y)]_i = \max[ l_i, \min[ y_i, u_i]]
  \eeqn
  and
  \beqn{cheap-projs2}
  [P_{\calF\cap\calB_k}(y)]_i = \max[l_i, x_{k,i}-\Delta_{k,i},
      \min[ y_i, x_{k,i}+\Delta_{k,i},u_i]].
  \eeqn
Note that \req{sL-def} can be rewritten as
\[
s_k^L = \mathcal{P}_{\mathcal{F}\cap\mathcal{B}_k} (x_k-g_k )- x_k 
= \mathcal{P}_{\mathcal{B}_k}(\mathcal{P}_{\mathcal{F}}(x_k-g_k ))-x_k
=\mathcal{P}_{\mathcal{B}_k}(x_k + d_k )-x_k. 
\]
\item
The uniform bound specified by AS.3 below is the only restriction made on
$B_k$ (beyond symmetry). This allows for a wide range of deterministic or random
approximations, such as (sampled) Barzilai-Borwein, safeguarded (limited-memory) 
quasi-Newton, random sketching or finite-difference approximations 
of $(s_k^L)^T B_k s_k^L$ using one evaluation of the gradient oracle.
Obviously, exact Hessians may also be used when available.
\item
The shifted quadratic model of the objective function given by
$g_k^Ts_k+\half s_k^TB s_k$ is minimized by $x_k+s_k^Q$ along the 
intersection of the first-order
direction $s_k^L$ with the trust region  $\calB_k$. In the vocabulary of
trust-region methods, it can therefore be interpreted as the ``Cauchy
point'' at iteration $k$ (see \cite[Sections~6.3 and 12.2.1]{ConnGoulToin00}). 
Since $s_k=s_k^Q$ satisfies \req{inside-and-decr}, an improved second-order 
step $s_k$ 
minimizing the quadratic model $g_k^Ts+\half s^TB s$ beyond $x_k+s_k^Q$ is possible 
and often beneficial. This is standard and well-tested practice in truncated Newton or second-order 
trust-region or adaptive-regularization algorithms.  The idea is to apply a conjugate-gradient or Lanczos method, which performs 
successive minimizations of the model in nested Krylov subspaces, until a desired 
accuracy threshold is reached (see for instance \cite{Toin81b,Stei83a} or 
\cite[Section~7.5]{ConnGoulToin00} for detailed descriptions and algorithms). Note that these 
techniques do not require computing a possibly approximate Hessian matrix, but only require 
its product with vectors.  
However, this remains optional in \algname\ (as is the evaluation of a nonzero $B_k$ and 
the computation of $s_k^Q$ itself) and the choice $B_k=0$ is possible when access to second-
order information is too expensive. Should this choice be made, \req{gamma-def} 
gives that $s_k^Q=s_k^L$ and the choice $s_k = s_k^L$ 
is always acceptable for \req{inside-and-decr}, in which case \algname\ is a purely first-order algorithm.
\item 
If the choice $B_k=0$ is made at all iterations and the problem is unconstrained in that 
$l_i=-\infty$ and $u_i=+\infty$ for all $i$, 
then $d_k=-g_k$ and thus \algname\ is nothing 
but the standard \tal{AdaGrad} algorithm. In short,
\[
\tim{ \algname\ = \tal{AdaGrad} + 2nd order + bound constraints,}
\]
the last two items being of course optional.
As a consequence, the stochastic complexity theory 
described below applies to \tal{AdaGrad} without any modification.

\end{enumerate}

\section{Convergence analysis}\label{section:convergence}

The convergence theory we are about to describe is based on the following 
assumptions.\\

\noindent
{\bf AS.1:} \textit{The function $f$ is twice continuously differentiable and the feasible region $\calF$ is not empty.}\\

\noindent
{\bf AS.2:} \textit{There exists a constant $L \geq0$ such that for all $x,y \in \Re^n$
\[
\|\nabla_x^1f(x)-\nabla_x^1f(y)\| \leq L\|x-y\|.
\]
}

\noindent
We stress that, although the existence of $L$ is assumed, the knowledge of 
its value is \textit{not} needed to run the algorithm.\\

\noindent
{\bf AS.3:} \textit{There exists a constant $\kappa_B\geq 1$ such that
$\|B_k\| \leq \kappa_B$ }for all $k \geq 0$. \\

\noindent
{\bf AS.4:} \textit{The objective function is bounded below on the feasible
domain, that is there exists a constant $\flow < f(x_0)$ such that $f(x)\geq
\flow$ for every $x \in \calF$.}\\

\noindent
{\bf AS.5:}
\textit{There exists a constant $\kappa_{Gg}>0$ such that
\[
\Econd{k}{|(G_k - g_k)^Ts_k|}
\leq \kappa_{Gg}^2 \Econd{k}{\|s_k\|^2}
\]
for all $k\geq0$.}\\

\noindent
This condition (which we have called "new1" in Table~\ref{table:adag-theory}) can be interpreted as a ``root mean square'' condition along the direction $s_k$ (see the comments after Theorem~\ref{theorem:true-convergence} below).
As far as the authors are aware, this condition is necessary if a Cauchy 
point is introduced to take second-order information into account because 
$\gamma_k$ then becomes a random variable.

Assumption AS.5 requires increasing accuracy when the iterates converge. This is needed in our analysis to estimate the conditional expectation of $G_k^Ts_k$  when the step is no longer a scaled multiple of the negative gradient, but can be quite general provided it produces a sufficient decrease in the quadratic model decrease.  Thus AS.5 is admittedly stronger than the typical bounded variance condition assumed in the simpler unconstrained case without second-order information, but it provides, in the significantly more challenging context studied in our paper, an optimal iteration complexity of $\mathcal{O}(\epsilon^{-2})$ rather than the $\mathcal{O}(\epsilon^{-4})$ obtained in the simpler context only assuming an unbiased oracle with bounded variance.

Assumption AS.5 is weaker than assuming 
\[
\Econd{k}{\|(G_k - g_k)\|^2}
\leq \kappa_{Gg}^4 \Econd{k}{\|s_k\|^2},
\]
as under this condition we have  
$$
\Econd{k}{|(G_k - g_k)^Ts_k|}\le \Econd{k}{\|G_k - g_k\|\|s_k\|}
\leq \sqrt{\Econd{k}{\|G_k - g_k\|^2}} \sqrt{\Econd{k}{\|s_k\|^2}}
\le \kappa_{Gg}^2 \Econd{k}{\|s_k\|^2}.
$$
In the paper \cite{cgt2012} it is shown  as in the deterministic case the condition 
$$
{\|(G_k - g_k)\|^2}
\leq \kappa_{Gg}^4{\|s_k\|^2},
$$
may be ensured by finite differences approximation of the gradient. 

The following ``linear descent'' lemma is a variant of
\cite[Lemma~2.1]{GratJeraToin24a} (which, in the unconstrained deterministic
context, uses a different optimality measure and a different
definition of $s_k^L$). 
It is also possible to replace AS.5 with the  condition 
\[
\Econd{k}{|(G_k - g_k)^Ts_k|}
\leq \kappa_{Gg}^2 \Econd{k}{\|s_k\|^2}+Q_k, \
\]
the rate of convergence then depends on the rate of convergence of $\{Q_k\}$ to zero (see the Appendix).

\llem{lemma:single-level-decrease}{
Suppose that AS.3 and AS.5 hold. Then, for $j\geq 0$,
\beqn{gen-decr}
\Econd{j}{G_j^Ts_j}
\leq
-\frac{\tau\varsigma^2}{2\kappa_B}\,\Econd{j}{d_j^T\Delta_j}
+ \kappa_s^2 \left(\frac 1 2 \kappa_B+\kappa_{Gg}^2\right)\,\Econd{j}{\|\Delta_j\|^2},
\eeqn
where $\Delta_j\in\Re^n$ is the vector whose $i$-th component is $\Delta_{j,i}$.
}

\proof{Consider any component $i\in\ii{n}$. We first note that \req{dk-def} 
and the contractive nature of the projection $P_\calF$ ensure 
that $|g_{j,i}|\geq |d_{j,i}|$.  Moreover,
\req{dk-def}, \req{wDk-def}, \req{Bk-def} and \req{sL-def} implies  
that either $|s_{j,i}^L| = |d_{j,i}|/w_{j,i}$
or $|s_{j,i}^L| = |d_{j,i}|$. In the first case, we may 
deduce that 
\beqn{L3.1:e1}
|g_{j,i}s_{j,i}^L|
\geq \frac{d_{j,i}^2}{w_{j,i}}
= w_{j,i} \frac{d_{j,i}^2}{w_{j,i}^2}
\geq \varsigma (s_{j,i}^L)^2,
\eeqn
while, in the latter case,
\beqn{L3.1:e2}
|g_{j,i}s_{j,i}^L|
\geq d_{j,i}^2 = (s_{j,i}^L)^2
= w_{j,i} \frac{d_{j,i}^2}{w_{j,i}}
\geq \varsigma \frac{d_{j,i}^2}{w_{j,i}}.
\eeqn
Combining the two cases and remembering that $\varsigma \in (0,1]$, we
see that, for $i\in\ii{n}$,
\beqn{L3.1:e3}
|g_{j,i}s_{j,i}^L| \geq \varsigma \frac{d_{j,i}^2}{w_{j,i}}
\tim{ and }
|g_{j,i}s_{j,i}^L| \geq \varsigma (s_{j,i}^L)^2.
\eeqn
Summing over all components $i\in\ii{n}$ and using the fact
that, by construction, $g_{j,i}s_{j,i}^L <0$ for all
$i$ then gives that 
\beqn{gsL}
g_j^Ts_j^L \leq - \varsigma \sum_{i=1}^n\frac{d_{j,i}^2}{w_{j,i}}
\tim{ and }
|g_j^Ts_j^L| \geq \varsigma \|s_j^L\|^2.
\eeqn
We now consider the quadratic model and suppose first 
that $\gamma_j < 1 $. Thus $(s_j^L)^TB_js_j^L>0$.  
We then deduce from \req{gamma-def}, AS.3 and \req{gsL} that 
\beqn{DL-4}
g_j^Ts_j^Q+\half (s_j^Q)^TB_js_j^Q
= -\bigfrac{(g_j^Ts_j^L)^2}{2(s_j^L)^TB_js_j^L}
\leq - \bigfrac{\varsigma}{2\kappa_B}|g_j^Ts_j^L|
\leq -\bigfrac{\varsigma^2}{2\kappa_B}\sum_{i=1}^n\frac{d_{j,i}^2}{w_{j,i}}.
\eeqn
If now $\gamma_j=1$, then \req{gamma-def}
and \req{gsL} give that
\[
g_j^Ts_j^Q+\frac 1 2  (s_j^Q)^TB_js_j^Q
= g_j^Ts_j^L+\frac 1 2 (s_j^L)^TB_js_j^L
\le \frac 1 2 g_j^Ts_j^L 
< - \frac{\varsigma}{2} \sum_{i=1}^n\frac{d_{j,i}^2}{w_{j,i}}
\]
and \req{DL-4} then again follows from the bounds
$\kappa_B \ge 1$ and $\varsigma\leq 1$. Thus, successively using 
the third part of \req{inside-and-decr}, AS.3, the second part of \req{inside-and-decr} and \req{DL-4},
\begin{align*}
g_j^Ts_j
&  = g_j^Ts_j + \half s_j^TB_js_j - \frac 1 2 s_j^TB_js_j\\*[1.5ex]
& \le \tau\left(g_j^Ts_j^Q + \frac{1}{2}(s_j^Q)^TB_js_j^Q\right)
   + \frac{1}{2}\kappa_B \|s_j\|^2 \\
& \le \tau\left(g_j^Ts_j^Q + \frac{1}{2}(s_j^Q)^TB_js_j^Q\right)
    +\frac{1}{2} \kappa_s^2\kappa_B \|\Delta_j\|^2\\
& \le -\frac{\tau\varsigma^2}{2\kappa_B} \bigsum_{i=1}^n\frac{d_{j,i}^2}{w_{j,i}}
       + \frac{1}{2} \kappa_s^2 \kappa_B \|\Delta_j\|^2,
\end{align*}
which, with the second part of \req{wDk-def}, gives that
\beqn{gen-decr0}
g_j^Ts_j 
\leq -\frac{\tau\varsigma^2}{2\kappa_B}d_j^T\Delta_j
     + \frac 1 2  \kappa_s^2 \kappa_B \|\Delta_j\|^2.
\eeqn
But, using AS.5 and the second part of \req{inside-and-decr},
\[
\begin{array}{lcl}
\Econd{j}{G_j^T s_j}
& =& \Econd{j}{g_j^Ts_j}+\Econd{j}{(G_j-g_j)^Ts_j}\\*[1ex]
&\leq& \Econd{j}{g_j^Ts_j}+|\Econd{j}{(G_j-g_j)^Ts_j}|\\*[1ex]
&\leq& \Econd{j}{g_j^Ts_j}+\kappa_{Gg}^2|\Econd{j}{\|s_j\|^2}|\\*[1ex]
&\leq& \Econd{j}{g_j^Ts_j}+\kappa_{Gg}^2\kappa_s^2\Econd{j}{\|\Delta_j\|^2}\\
\end{array}
\]
and \req{gen-decr} follows from \req{gen-decr0}.% and \req{sL-def}. k
} % epr

\noindent
Lemma~\ref{lemma:single-level-decrease} is crucial for the proof of our 
complexity bounds, which also depends on two technical results.  

\llem{tech-tan}{
  Suppose that, for $t>0$ and $a,c>0$,
  \beqn{abc}
  at \le b + c\,\log(t).
  \eeqn
  Then, 
  \beqn{first-bound}
  t \le \frac{2b}{a}+ \frac{2c}{a}\left[\log\left(\frac{2c}{a}\right)-1\right].
  \eeqn
  }

  \proof{The concavity of the logarithm at $2$ gives that, for all $x>0$
  \[
  \log(x)\le \log(2)+\frac{x-2}{2} =\frac{x}{2}+\log(2)-1.
  \]
  Applying this inequality for $x = at/c$ then gives that
  \private{%%%%
  \[
  \log\left(\frac{at}{c}\right)\le \frac{at}{2c} + \log(2) -1
  \]
  \[
  \log(t) 
  \le \frac{at}{2c} + \log\left(\frac{c}{a}\right) + \log(2) -1
  \]
  }%%%% end private
  \[
  \log(t)=\frac{at}{2c} +\log\left(\frac{2c}{a}\right)-1,
  \]
  which we may substitute  in the inequality $2at\le 2b+2c\log(t)$ obtained by multiplying \req{abc} by $2$, yielding that
  \[
  2t
  \le \frac{2b}{a}+\frac{2c}{a}\left[\frac{at}{2c} + \log\left(\frac{2c}{a}\right)-1\right]
   =\frac{2b}{a} + t + \frac{2c}{a}\left[\log\left(\frac{2c}{a}\right)-1\right],
  \]
  and \req{first-bound} follows by substracting $t$ on both sides.
}%epr

\llem{magical}{
  Let $\{a_j\}_{j\ge 0}$ be a sequence of non-negative numbers and let
  $b_j=\sum_{i=0}^j a_i$.  Then
  \[
  \sum_{j=0}^k \frac{a_j}{\varsigma+b_j} \le
  \log\left(1+\frac{1}{\varsigma}b_k\right).
  \]
}

\proof{See \cite{ward2020adagrad}. }

\noindent
Our first complexity result now consider the global rate of
convergence of $\|d_k\|$ to zero. 
%It is expressed using the 
%expectation conditioned by $\calA$, which we denote by the symbol $\EA{\cdot}$.

\lthm{theorem:convergence}{
Suppose that AS.1--AS.5 hold and that the
\algname\ algorithm is applied to problem \req{problem}. Then
\beqn{gradbound}
\E{\average_{j\in\iiz{k}} \|d_j\|} \le \frac{\kap{conv}}{\sqrt{k+1}},
\eeqn
where 
\beqn{kapconv-def}
\kap{conv} =\frac{4\kappa_B}{\tau\varsigma^2}\left\{1+\Gamma_0+2n\kappa_*
\left[\log\left(\frac{8n\kappa_*\kappa_B}{\tau\varsigma^{5/2}}\right)-1\right]\right\}
%\kap{conv}=\sqrt{\varsigma}\,
%\max\left[e^{\frac{1 + f(x_0)-\flow}{2n\kappa_*}},
%\frac{64n^2\kappa_*^2\kappa_B^2}{\tau^2\varsigma^5}\right]
\tim{ with } 
\kappa_* = \kappa_s^2(\kappa_{Gg}^2+\half (\kappa_B + L)).
\eeqn
}

\proof{
Consider an iteration $j$ and first note that
\[
\|s_j\|^2
%\leq \kappa_s^2\|s_j^L\|^2
\leq \kappa_s^2\|\Delta_j\|^2
= \kappa_s^2\bigsum_{i=1}^n\bigfrac{d_{j,i}^2}{w_{j,i}^2}
\]
because of the second part of \req{inside-and-decr} %\req{sL-def}, \req{Bk-def} 
and \req{wDk-def}. Lemma~\ref{lemma:single-level-decrease},
AS.1 and AS.2 then give that
\beqn{conv1-12}
\begin{array}{lcl}
\Econd{j}{f(x_{j+1})}
& \leq & f(x_j) + \Econd{j}{G_j^Ts_j} + \bigfrac{L}{2}\Econd{j}{\|s_j\|^2} \\*[1ex]
& \leq & f(x_j) -\bigfrac{\tau\varsigma^2}{2\kappa_B}\Econd{j}{\bigsum_{i=1}^n\bigfrac{d_{j,i}^2}{w_{j,i}}}
+\kappa_s^2(\kappa_{Gg}^2+ \frac{1}{2}\kappa_B + \frac 1 2 L)\,\Econd{j}{\bigsum_{i=1}^n\bigfrac{d_{j,i}^2}{w_{j,i}^2}}.
\end{array}
\eeqn
Defining $\kappa_* =  \kappa_s^2(\kappa_{Gg}^2+\half \kappa_B + \half L)$
and using the law of total expectation gives that
\[
\E{f(x_{j+1})}\leq \E{f(x_j)} - \bigfrac{\tau\varsigma^2}{2\kappa_B}\E{\bigsum_{i=1}^n\bigfrac{d_{j,i}^2}{w_{j,i}}}
+ \kappa_*\E{\bigsum_{i=1}^n\bigfrac{d_{j,i}^2}{w_{j,i}^2}},
\]
and therefore, summing for $j\in\iiz{k}$ for $k$ fixed, that
\beqn{conv-decr2}
\bigfrac{\tau\varsigma^2}{2\kappa_B}\E{\bigsum_{j=0}^k\bigsum_{i=1}^n\bigfrac{d_{j,i}^2}{w_{j,i}}}
\leq f(x_0)-\flow + \kappa_*\E{\bigsum_{j=0}^k\bigsum_{i=1}^n\bigfrac{d_{j,i}^2}{w_{j,i}^2}}.
\eeqn
Observe now that, if $p= \argmax_{q\in\ii{n}}w_{j,q}$,
\[
w_{j,i} 
\le {\bigmax_{q\in\ii{n}}w_{j,q}= w_{j,p}}=\sqrt{\varsigma + \sum_{\ell=0}^j d_{\ell,p}^2}
\le \sqrt{\varsigma + \sum_{\ell=0}^j\|d_\ell\|^2}
\]
for all $i\in\ii{n}$ and thus that
\beqn{anotherbound}
\bigsum_{j=0}^k\frac{\|d_j\|^2}{\sqrt{\varsigma+\sum_{\ell=0}^j\|d_\ell\|^2}}
{\le\bigsum_{j=0}^k\bigsum_{i=1}^n\bigfrac{d_{j,i}^2}{\bigmax_{q\in\ii{n}}w_{j,q}}}
\le \bigsum_{j=0}^k\bigsum_{i=1}^n\bigfrac{d_{j,i}^2}{w_{j,i}}.
\eeqn
We may now apply Lemma~\ref{magical} to
derive that, for each $i\in\ii{n}$,
\[
\begin{array}{lcl}
\E{\bigsum_{i=1}^n\bigsum_{j=0}^{k}\frac{d_{j,i}^2}{w_{j,i}^2}}
& = & \displaystyle \E{\sum_{i=1}^n\sum_{j=0}^{k}\frac{d_{j,i}^2}{\varsigma+   \sum_{\ell=0}^j d_{j,\ell}^2}}\\
& \leq & \displaystyle n\E{\max_{i\in\ii{n}}\log\left(1+\frac{1}{\varsigma}\sum_{j=0}^{k}d_{j,i}^2\right)}\\
& \leq & \displaystyle n\E{\log\left(1+\frac{1}{\varsigma}\sum_{i=1}^n\sum_{j=0}^{k}d_{j,i}^2\right)},
\end{array}
\]
and hence, using AS.4 to define $\Gamma_0 = f(x_0)-\flow>0$,
\req{conv-decr2}  and \req{anotherbound} gives that
\beqn{abound}
\bigfrac{\tau\varsigma^2}{2\kappa_B}\E{\bigsum_{j=0}^k\frac{\|d_j\|^2}{\sqrt{\varsigma+\sum_{\ell=0}^j\|d_\ell\|^2}}}
\leq \Gamma_0 + n\kappa_*\E{\log\left(1+\frac{1}{\varsigma}\sum_{j=0}^{k}\|d_j\|^2\right)}
\eeqn
in turn ensuring that
\beqn{abound2}
\begin{aligned}
\bigfrac{\tau\varsigma^2}{2\kappa_B}\E{\bigsum_{j=0}^k\frac{\|d_j\|^2}{\sqrt{\varsigma+\sum_{\ell=0}^k\|d_\ell\|^2}}}
&\le \Gamma_0 + n\kappa_*\E{\log\left(1+\frac{1}{\varsigma}\sum_{j=0}^{k}\|d_j\|^2\right)}\\
&= \Gamma_0+2n\kappa_*\E{\log\left(\sqrt{1+\frac{1}{\varsigma}\sum_{j=0}^{k}\|d_j\|^2}\right)}\\
&\leq \Gamma_0 +2n\kappa_*\log\left(\E{\sqrt{1+\frac{1}{\varsigma}\sum_{j=0}^{k}\|d_j\|^2}}\right),
\end{aligned}
\eeqn
where we used Jensen's inequality and the concavity of the logarithm to deduce the last inequality.
Now
\[
\bigsum_{j=0}^k\frac{\|d_j\|^2}{\sqrt{\varsigma+\sum_{\ell=0}^k\|d_\ell\|^2}}
\private{
= \frac{\varsigma+\sum_{j=0}^k\|d_j\|^2}  
{\sqrt{\varsigma+\sum_{\ell=0}^k\|d_\ell\|^2}} 
-\frac{\varsigma}{\sqrt{\varsigma+\sum_{\ell=0}^k\|d_\ell\|^2}}   }
= \sqrt{\varsigma+\sum_{\ell=0}^k\|d_\ell\|^2}
-\frac{\varsigma}{\sqrt{\varsigma+\sum_{\ell=0}^k\|d_\ell\|^2}}
\]
and thus
\[
\sqrt{\varsigma}\;\E{\sqrt{1+\frac{1}{\varsigma}\sum_{\ell=0}^k\|d_\ell\|^2}}
=\E{\sqrt{\varsigma+\sum_{\ell=0}^k\|d_\ell\|^2}}
\le \E{\bigsum_{j=0}^k\frac{\|d_j\|^2}{\sqrt{\varsigma+\sum_{\ell=0}^k\|d_\ell\|^2}}}
+\sqrt{\varsigma}.
\]
Substituting \req{abound2} in this inequality and using that $\tau\varsigma^2\sqrt{\varsigma}\le 1 < 2\kappa_B$ then gives that
\[
\bigfrac{\tau\varsigma^2\sqrt{\varsigma}}{2\kappa_B}\E{\sqrt{1+\frac{1}{\varsigma}\sum_{\ell=0}^k\|d_\ell\|^2}}
\le 1+ \Gamma_0 + 2n\kappa_*\log\left(\E{\sqrt{1+\frac{1}{\varsigma}\sum_{\ell=0}^k\|d_\ell\|^2}}\right),
\]
which can be rewritten as
\[
a\,t_k \le b + c\,\log(t_k)
\]
with
\[
t_k = \E{\sqrt{1+\frac{1}{\varsigma}\sum_{\ell=0}^k\|d_\ell\|^2}}\ge1,
\ms
a = \bigfrac{\tau\varsigma^\sfrac{5}{2}}{2\kappa_B},
\ms
b = 1 + \Gamma_0 
\tim{ and }
c = 2n\kappa_*.
\]
We may then use Lemma~\ref{tech-tan} and deduce that
\[
t_k\le \frac{2b}{a}+\frac{2c}{a}\left[\log\left(\frac{2c}{a}\right)-1\right]
= \frac{4\kappa_B}{\tau\varsigma^{5/2}}\left\{1+\Gamma_0+2n\kappa_*
\left[\log\left(\frac{8n\kappa_*\kappa_B}{\tau\varsigma^{5/2}}\right)-1\right]\right\},
\]
which gives that
\[
\begin{aligned}
\E{\sqrt{\sum_{\ell=0}^k\|d_\ell\|^2}}
& \le \sqrt{\varsigma}\,\E{\sqrt{1+\frac{1}{\varsigma}\sum_{\ell=0}^k\|d_\ell\|^2}}\\
&\le \frac{4\kappa_B}{\tau\varsigma^2}\left\{1+\Gamma_0+2n\kappa_*
\left[\log\left(\frac{8n\kappa_*\kappa_B}{\tau\varsigma^{5/2}}\right)-1\right]\right\}\\
&= \kap{conv}.
\end{aligned}
\]
Dividing by $\sqrt{k+1}$ and using the inequality 
\[
\frac{1}{k+1}\sum_{j=0}^k\|d_j\| \leq \frac{1}{\sqrt{k+1}}\sqrt{\sum_{j=0}^k \|d_j\|^2}
\]
finally gives
\[
\E{\average_{j\in\iiz{k}}\|d_j\|}
\leq \E{\sqrt{\average_{j\in\iiz{k}}\|d_j\|^2}}
\leq  \frac{\kap{conv}}{\sqrt{k+1}},
\]
as requested.
}%epr

\noindent
In the deterministic case, Theorem~\ref{theorem:convergence} provides
a bound on the complexity of solving the bound-constrained problem
\req{problem} for which we consider the standard optimality measure
\cite[Section~12.1]{ConnGoulToin00}
$\|\Xi_k\|$, where
\[
\Xi_k \eqdef P_\calF(x_k-G_k)-x_k
\]
(if the problem is unconstrained, $\|\Xi_k\|= \|G_k\|$).

\lcor{theorem:convergence-deterministic}{
Suppose that AS.1--AS.4 hold, that the
\algname\ algorithm is applied to problem \req{problem} 
%starting from a non-critical $x_0$ with$\varsigma < \|d_0\|^2$ 
and that $g_j = G_j$ for all $j$. Then
\beqn{gradbound-deter}
\min_{j\in\iiz{k}} \|\Xi_j\|
\leq \average_{j\in\iiz{k}} \|\Xi_j\|
\leq \frac{\kap{conv}}{\sqrt{k+1}},
\eeqn
where the constant $\kap{conv}$ is computed as in
Theorem~\ref{theorem:convergence} using the value $\kappa_{Gg} = 0$.
}

\proof{
The result trivially follows from Theorem~\ref{theorem:convergence} with $\kappa_{Gg}=0$.
}

\section{First-order optimality and stochastic projected gradients with bounds}\label{section:spg}

\subsection{A more general framework}

Corollary~\ref{theorem:convergence-deterministic} gives the global rate of convergence of the \algname\ algorithm in the deterministic case where (obviously) $G_k = g_k$ and thus $\Xi_k=d_k$. Things are more complicated in the stochastic case because, as we prove in this section, this last inequality may no longer hold.  Thus, while Theorem~\ref{theorem:convergence} gives the rate of convergence of the \textit{approximate} criticality measure, obtaining a similar rate for the \textit{true} criticality measure is another question. Because our way to answer this question applies
not only to the \algname\ algorithm, but also to a wider class of stochastic projected-gradient methods, we now consider solving problem \req{problem} using the algorithmic framework given by \tal{StochProjGrad}.

\algo{spg-algo}{\tal{StochProjGrad}
    ($x_{\rm ini},l,u$)}{
\begin{description}
\item[Step 0: Initialization:]
Set  $x_0 = P_\calF(x_{\rm ini})$, $k=0$. 
\item[Step 1: Step:]
  Compute $g_k= g(x_k,\xi)$ a random approximation of $G_k=\nabla_x^1f(x_k)$, and a step $s_k$ such that
  $x_k+s_k \in \calF$,
\item[Step 4: Loop:] Set $x_{k+1} = x_k + s_k$, increment $k$ by one and go to Step~1.
\end{description} 
}
\vspace*{3mm}

\noindent
Moreover we will assume, in this subsection, that the expected approximate criticality measure 
\[
\Econd{k}{\|d_k\|} = \Econd{k}{\|P_\calF(x_k-g_k)-x_k\|}
\]
converges to zero (at least on average as in Theorem~\ref{theorem:convergence}). We are now interested in what can be deduced on $\E{\|\Xi_k\|} = \E{\|P_\calF(x_k-G_k)-x_k\|}$, the relevant criticality measure for problem \req{problem} in the stochastic case. Ideally, one would hope that the approximate gradient's distribution ensures coherence of the measures in the sense that
\beqn{the-condition}
\E{\|\Xi_k\|} \leq \kap{opt}\E{\|d_k\|}
\eeqn 
for a fixed $\kap{opt}>0$,
in which case $\E{\|\Xi_k\|}$ converges to zero at the same rate as $\E{\|d_k\|}$.
If we consider the unconstrained case ($\calF=\Re^n$) and assume that the gradient oracle is
unbiased (i.e. $\Econd{k}{g_k}=G_k$), then, using Jensen's inequality
and the convexity of the norm,
\[
\|\Xi_k\|  = \|G_k\| = \|\Econd{k}{g_k}\| \leq \Econd{k}{\|g_k\|}= \Econd{k}{\|d_k\|}.
\]
Taking the full expectation on both sides and using the law of total 
expectation then shows that \req{the-condition} always holds with $\kap{opt}=1$.
The situation is very different if bounds are present, even if the
gradient oracle is unbiased.  To see this, consider the following
one dimensional example, where $\calF = [0,+\infty)$ and
\[
x_k = \frac{1}{k+1}
\tim{ and }
g_k = \left\{\begin{array}{ll} 1 &\tim{with probability} p_k = \bigfrac{1}{k+1}+\bigfrac{1}{(k+1)^2}\\
0 &\tim{with probability} 1-p_k,
\end{array}\right.
\]
for $k> 0$. Also define $G_k = \Econd{k}{g_k} = p_k$ (so that the gradient
oracle is unbiased). One then easily verifies that 
$|G_{k+1}-G_k| \leq 3\,|x_{k+1}-x_k|$, so that AS.2 holds with $L=3$.
Indeed, since $G_k= x_k+x_k^2$, we have that
\[
|G_{k+1}-G_k|
\leq |x_{k+1}-x_k+x^2_{k+1}-x^2_k|\le |x_{k+1}-x_k+(x_{k+1}-x_k)(x_{k+1}+x_k)|\le   3\, |x_{k+1}-x_k|.
\]
These definitions give that 
\[
|\Xi_k| 
= \left|P_{[0,+\infty)}\left(\bigfrac{1}{k+1}-p_k\right)-\bigfrac{1}{k+1}\right|
= \left|P_{[0,+\infty)}\left(\bigfrac{-1}{(k+1)^2}\right)-\bigfrac{1}{k+1}\right|
= \bigfrac{1}{k+1}
\]
and
\[
\begin{array}{lcl}
\Econd{k}{|d_k|}
& = & \left|P_{[0,+\infty)}\left(\bigfrac{1}{k+1}-1\right)-\bigfrac{1}{k+1}\right|\,p_k
   + \left|P_{[0,+\infty)}\left(\bigfrac{1}{k+1}-0\right)-\bigfrac{1}{k+1}\right|\,(1-p_k)\\*[2ex]
& = & \bigfrac{1}{k+1}\,p_k + 0 \,(1-p_k)\\*[2ex]
& = & \bigfrac{1}{(k+1)^2} + \bigfrac{1}{(k+1)^3}.
\end{array}
\]
Applying now the law of total expectation, we see that
\[
\lim_{k\rightarrow\infty}\frac{\E{|d_k|}}{\E{|\Xi_k|}} = 0,
\]
preventing \req{the-condition} (a "coherently distributed" gradient oracle)
to hold for a fixed $\kap{opt}>0$.
We therefore conclude that, in general, 
\textit{the rate of decrease of $\E{\|d_k\|}$ does} not 
\textit{translate to a similar rate of decrease for $\E{\|\Xi_k\|}$, 
even for unbiased gradient oracles}.

Fortunately, the situation can be improved by strengthening the condition on the
gradient accuracy, even if the gradient oracle is biased. 
This is the object of the next lemma.

\llem{lemma:bias}{
For each $k\geq0$ and each $i\in\ii{n}$, we have that
\beqn{bias2}
\E{\|\Xi_k\|}
\leq \E{\|d_{k}\|} + \E{\|g_{k}-G_{k}\|}.
\eeqn
Moreover, if $\Econd{k}{\|g_{k}-G_{k}\|}\leq \kap{err} \,\Econd{k}{\|d_{k}\|}$ for some $\kap{err}\geq 0$,
then
\beqn{condbiased2}
\E{\|\Xi_k\|} \leq (1+\kap{err})\,\E{\|d_k\|}.
\eeqn
}

\proof{
Consider an arbitrary $i\in\ii{n}$ and note that
\beqn{Xik-def}
\Xi_{k,i} = P_{i}(x_{k,i}-G_{k,i})-x_{k,i}.
\eeqn
where $P_{i} = P_{[l_{i},u_{\ell_i}]}$. Since $P_{i}$ is a contractive map, we have that
\[
\begin{array}{l}
(P_{i}(x_{k,i}-G_{k,i})-x_{k,i})^2 - (P_{i}(x_{k,i}-g_{k,i})-x_{k,i})^2\\*[1ex]
\hspace*{20mm} = \left[ P_{i}(x_{k,i}-G_{k,i})-x_{k,i}+P_{i}(x_{k,i}-g_{k,i})-x_{k,i}\right] 
\times \\*[1ex]
\hspace*{30mm}\left[P_{i}(x_{k,i}-G_{k,i})-P_{i}(x_{k,i}-g_{k,i})\right ]\\*[1ex]
\hspace*{20mm} \leq \left| P_{i}(x_{k,i}-G_{k,i})-x_{k,i}+P_{i}(x_{k,i}-g_{k,i})-x_{k,i}\right|\,|g_{k,i}-G_{k,i}|.
\end{array}
\]
Now, again using the contractivity of $P_{i}$,
\[
\begin{array}{l}
\left| P_{i}(x_{k,i}-G_{k,i})-x_{k,i}+P_{i}(x_{k,i}-g_{k,i})-x_{k,i}\right|\\*[1ex]
\hspace*{20mm} = \big| 2(P_{i}(x_{k,i}-g_{k,i})-x_{k,i})\\*[1ex]
\hspace*{30mm}+ (P_{i}(x_{k,i}-G_{k,i})-x_{k,i})-(P_{i}(x_{k,i}-g_{k,i})-x_{k,i})\big|\\*[1ex]
\hspace*{20mm} \leq 2\left|P_{i}(x_{k,i}-g_{k,i})-x_{k,i}\right|+\left|P_{i}(x_{k,i}-G_{k,i})-P_{i}(x_{k,i}-g_{k,i})\right|\\*[1ex]
\hspace*{20mm} \leq 2\left|P_{i}(x_{k,i}-g_{k,i})-x_{k,i}\right|+\left|g_{k,i}-G_{k,i}\right|\\
\end{array}
\]
and thus
\[
\begin{array}{l}
(P_{i}(x_{k,i}-G_{k,i})-x_{k,i})^2 - (P_{i}(x_{k,i}-g_{k,i})-x_{k,i})^2\\*[1ex]
\hspace*{20mm} \leq 2\left|P_{i}(x_{k,i}-g_{k,i})-x_{k,i}\right|\,\left|g_{k,i}-G_{k,i}\right|+\left|g_{k,i}-G_{k,i}\right|^2,
\end{array}
\]
which, using \req{dk-def} and \req{Xik-def}, gives that
\[
\Xi_{k,i}^2 
\leq d_{k,i}^2 + 2|d_{k,i}|\left|g_{k,i}-G_{k,i}\right| + \left|g_{k,i}-G_{k,i}\right|^2
\]
and therefore, summing for $i\in\ii{n}$ and using the Cauchy-Schwartz inequality, that
\[
\begin{array}{lcl}
\|\Xi_k\|^2 
& = & \|d_k\|^2 + 2|d_k|^T|g_k-G_k| + \|g_k-G_k\|^2\\*[2ex]
& \leq & \|d_k\|^2 + 2\|d_k\|\,\|g_k-G_k\| + \|g_k-G_k\|^2\\*[1.5ex]
& = & \Big( \|d_k\|+\|g_k-G_k\| \Big)^2.
\end{array}
\]
This yields that
\beqn{aproperty2}
\|\Xi_k\|\leq \|d_k\| + \|g_k-G_k\|.
\eeqn
and \req{bias2} follows by taking conditional expectations on both sides.
Moreover, 
if $\Econd{k}{\|g_{k}-G_{k}\|}\leq \kap{err} \,\Econd{k}{\|d_{k}\|}$, \req{condbiased2} 
results from \req{aproperty2}
\private{
\[
\Econd{k}{\|\Xi_k\|} \leq \Econd{k}{\|d_k\|} + \Econd{k}{\|g_{k}-G_{k}\|}
\leq (1+\kap{err})\Econd{k}{\|d_k\|}
\]
}%
and the tower property.
}

\noindent
Thus we see from \req{bias2} that the convergence of $\E{\|\Xi_k\|}$ to zero is still guaranteed if $\E{\|g_k-G_k\|}$ converges to zero along with $\E{\|d_k\|}$.  Moreover, if it does so at the same speed, the rate of decrease of $\E{\|d_k\|}$ dominates, as shown by \req{condbiased2}.

This result also indicates what can happen if nothing is known on the error in the gradient.  The optimality measure of the approximate problem goes to zero (as in Theorem~\ref{theorem:convergence}) but the 
true measure could remain at a level given by the second term of \req{bias2}.
To interpret this term, observe that, using Jensen's inequality and the concavity of the square root, 
\beqn{STD}
\Econd{k}{\|g_k-G_k\|}
= \Econd{k}{\sqrt{\|g_k-G_k\|^2}}
\leq \sqrt{\Econd{k}{\|g_k-G_k\|^2}}
= \RMSE_k
\eeqn
for all $k\geq 0$, where $\RMSE_k$ is the conditional
root mean square error (RMSE) of the gradient oracle at iteration $k$. Thus, by the law of total expectation,
$\E{\|g_k-G_k\|}\leq \E{\RMSE_k}$, and $\limsup_{k\rightarrow\infty} \E{\RMSE_k}$ gives an asymptotic upper bound on the criticality default (the deviation of the criticality measure from zero).

\subsection{Application to the \tal{ADAGB2} algorithm}

The discussion of the previous subsection finally allows us to rephrase 
Theorem~\ref{theorem:convergence} to cover convergence of the \algname\ algorithm on problem \req{problem} 
in three progressively more general scenarii (condition \req{bias-condition} 
below was called "new2" in Table~\ref{table:adag-theory}).

\lthm{theorem:true-convergence}{
Suppose that AS.1--AS.5 hold and that the
\algname\ algorithm is applied to problem \req{problem}. 
Then
\beqn{gradbound-true}
\average_{j\in\iiz{k}} \E{\|\Xi_j\|}
\le \frac{\kap{conv}\kap{1}}{\sqrt{k+1}}+\kappa_2\average_{j\in\iiz{k}}\E{\|g_k-G_k\|} +\kap{3}
\eeqn
where $\kap{conv}$ is defined in \req{kapconv-def}
 and
\begin{description}
\item[coherently distributed:] {$\kappa_1 = \kap{opt}$} and $\kappa_2= \kap{3}= 0$ if
  \[
  \E{\|\Xi_k\|} \leq \kap{opt}\E{\|d_k\|}
  \]
  for some constant $\kap{opt}>0$ and all $k\geq 0$;
\item[controlled error:] $\kappa_1 = 1+\kap{err1}$ and $\kappa_2= \kap{3}= 0$ if 
\beqn{bias-condition}
\Econd{k}{\|g_k-G_k\|}\leq \kap{err1} \,\Econd{k}{\|d_k\|}
\eeqn
for some constant $\kap{err1}>0$ and all $k\geq 0$;
\item[probabilistic error control:]
$\kappa_1 = (1+\kap{pec})$, $\kap{2}=0$ and $\kap{3} = (1-p)\kap{err2}$ 
if
\begin{eqnarray}
& &  \big\|g_k-G_k\big\|<\kap{err2}  \label{prob-condition1}\\
& & \mathbb{P}_k\big[{\cal{E}}_k\big]\ge p  \mbox{ with event } {\cal{E}}_k=\{\big\|g_k-G_k\big\|
\le \kap{pec} \|d_k\|\} \label{prob-condition2}
\end{eqnarray}
for some $p\in (0,1)$, $\kap{err2},\kap{pec}\ge0$ and all $k\geq 0$;
\item[general:] $\kappa_1 = \kappa_2 = 1$, \new{$\kappa_3=0$}, otherwise.
\end{description}
}

\proof{
The results for the ``coherently distributed'',  ``controlled error'' and ``general'' scenarii are obtained by combining Theorem~\ref{theorem:convergence} with Lemma~\ref{lemma:bias} and (for the ``coherently distributed'' and ``controlled error'' scenarii) from the comments following it.

For analyzing the ``probabilistic error control'' scenario, 
first note that, using the contractive nature of projection,
\[\big\|\Xi_k\big\|=\|P_\calF(x_k-g_k)-x_k +
 P_\calF(x_k-G_k)-P_\calF(x_k-g_k)\| 
 \le \|d_k\| + \|g_k-G_k\|.
 \]
If ${\mathbbm{1}}({\cal{E}}_k)=1$ then this inequality implies that
\[
\big\|\Xi_k\big\| \le (\kap{pec}+1) \|d_k\|,
\]
and it also implies that in all cases,
$$\big\|\Xi_k\big\|
\le \|d_k\|+\kap{err2}  
\le (\kap{pec}+1)\|d_k\|+ \kap{err2}.
$$ 
Thus, letting $\overline{ {\cal{E}}_k}=\{\big\|g_k-G_k\big\|
> \kap{pec} \|d_k\|\}$ be the complement of $\cal{E}$, we deduce that
\begin{eqnarray*}
\mathbb{E}_k[\|\Xi_k\|]&=& P_k({\cal{E}}_k)\, \mathbb{E}_k[\|\Xi_k\| \, |\, {\cal{E}}_k]+ P_k(\overline{ {\cal{E}}_k})\, \mathbb{E}_k[\|\Xi_k\|| \,\overline{ {\cal{E}}_k}]\\
& \le & (\kap{pec}+1)\mathbb{E}_k[\|d_k\|] +(1-p) \kap{err2}.
\end{eqnarray*}
Using now the tower property and (\ref{gradbound})
we obtain that
$$\average_{j\in\iiz{k}} \mathbb{E}{\|\Xi_j\|}
\le \frac{\kap{conv}(\kap{pec}+1)}{\sqrt{k+1}}+(1-p)\kap{err2}.
$$
}

If the function $f$ is subject to independent and identically distributed noise, then gradient can be estimated 
by standard sample averaging approximation techniques and 
the probabilistic condition in (\ref{prob-condition2}) can be enforced by  the Bernstein inequality \cite{BeGu22, MT20}.
The term depending on inaccurate approximations of the gradient is damped by the factor $(1-p)$ and its magnitude can be made arbitrarily small by increasing $p$.

Theorem~\ref{theorem:true-convergence} gives the desired fast convergence to zero of the
first-order optimality measure associated with problem \req{problem}
in the coherently distributed and controlled error cases.  In particular, if
the problem is unconstrained and the gradient oracle is unbiased,
then the coherently distributed case applies with $\kap{opt}=1$ and we
\textit{recover the $\calO(\epsilon^{-2})$ complexity bound for}
\tal{AdaGrad} obtained in the references mentioned in the
introduction \textit{assuming the condition given by
AS.5.}  The same is true in the bound-constrained case in 
the controlled-error scenario.

Theorem~\ref{theorem:true-convergence} and the discussion of 
the previous section also indicate that the 
true measure could remain at a level given by $\limsup_{k\rightarrow\infty}\beta_k$
where $\beta_k =\average_{j\in\iiz{k}}\E{\|g_k-G_k\|}$
is bounded above by the average $\RMSE_k$ of 
the gradient oracle, the average being taken on the first $k$ iterations.
If $\beta_k$ is bounded for large $k$ by a
(hopefully) small positive constant, then the optimality measure is only constrained 
to fall below this constant. However, convergence to true optimality must still 
happen if the sequence $\{\beta_k\}$ converges to zero\footnote{Which does not 
necessarily requires the convergence of the sequence $\{\|g_k-G_k\|\}$ to zero.}, 
albeit possibly at a rate slower than $\calO(1/\sqrt{k+1})$.

Moreover, reformulating the "new2" condition \req{bias-condition} using 
the bound \req{STD} allows us to state the following direct consequence of 
Theorem~\ref{theorem:true-convergence}.

\lcor{corollary:RMSE}{
Suppose that AS.1--AS.5 hold and that the
\algname\ algorithm is applied to problem \req{problem}. Suppose also that there 
exists a constant $\kap{err}\geq 0$ such that, for all $k\geq 0$,
\[
\RMSE_k \leq \kap{err} \,\Econd{k}{\|d_k\|}
\]
where $\RMSE_k$ is the conditional root mean square error defined by the last equality of \req{STD}.
Then
\[
\average_{j\in\iiz{k}} \E{\|\Xi_j\|} \leq \frac{(1+\kap{err})\kap{conv}}{\sqrt{k+1}},
\]
where $\kap{conv}$ defined by \req{kapconv-def}.
}

We finally state a complexity result in probability simply
derived from Theorem~\ref{theorem:true-convergence} for its first two scenarii.

\lcor{prob-complexity}{
Suppose that the conditions of Theorem~\ref{theorem:true-convergence}
hold for the coherently distributed or controlled-error cases.  Then,
given $\delta\in (0,1)$, one has that 
\[
\Prob{\min_{j\in\iiz{k}} \|\Xi_j\| \leq \epsilon} \geq 1 -\delta
\tim{ for }
k > \left(\frac{\kap{1}\kap{conv}}{\delta\epsilon}\right)^2,
\]
where $\kap{conv}$ is given by \req{kapconv-def}.}

\proof{
 Using Markov's inequality
 and \req{gradbound-true} with $\kap{2}=\kap{3}=0$, we have that
\[
\Prob{\min_{j\in\iiz{k}} \|\Xi_j\| \leq \epsilon}
\geq \Prob{\frac{1}{k}\bigsum_{j=0}^k\|\Xi_j\|\leq \epsilon}
\geq 1 -
\frac{1}{\epsilon}\,\E{\frac{1}{k}\bigsum_{j=0}^k\|\Xi_j\|}
  \geq 1-\frac{\kap{1}\kap{conv}}{\epsilon\sqrt{k+1}}.
\]
The desired conclusion follows.
}

\section{Conclusions and perspectives}\label{section:conclusion}

We have introduced the \algname\ algorithm for bound-constrained 
stochastic optimization which is also capable of using available 
second-order information. When second-order information is not 
used and the problem is unconstrained, \algname\ subsumes the 
\tal{AdaGrad} algorithm. 

We have shown that, given $\delta\in(0,1)$ and $\epsilon\in(0,1]$, 
\textit{the \algname\ algorithm needs at most $\calO(\epsilon^{-2})$ 
iterations to ensure an $\epsilon$-approximate first-order 
critical point of the bound-constrained problem \req{problem} 
with probability at least $1-\delta$}, provided that Assumption
AS.5 holds  and the average 
RMSE $\beta_k$ is sufficiently small. Should this condition fail, 
we have shown that the optimality default is bounded above for 
large $k$ by the limit superior of the average oracle's RMSE.

We have also discussed the relation between the standard optimality 
measure for bound-constrained induced by the projected gradient in 
the stochastic case, and have shown that, in general, the unbiased 
nature of the gradient oracle is not sufficient to ensure convergence 
on the true constrained problem at the optimal rate, even if such a 
convergence occurs for the approximate one.

Some interesting theoretical questions remain open at this point, 
including the extension of the results to $(L_0,L_1)$-smooth functions, 
the use of momentum, ensuring second-order optimality and the handling of more general constraints.

{\footnotesize
\section*{\footnotesize Acknowledgements}
Serge Gratton and Philippe Toint are grateful to Defeng Sun, Xiaojun
Chen and Zaikun Zhang of the Department of Applied Mathematics of
the Hong-Kong Polytechnic University for their support during a research
visit in the fall 2024. Philippe Toint also acknowledges the support
of DIEF (UNIFI) for a visit in October 2024. The research  of Stefania 
Bellavia and Benedetta Morini  was partially supported by INDAM-GNCS 
through Progetti di Ricerca 2023 and by PNRR - Missione 4 Istruzione 
e Ricerca - Componente C2 Investimento 1.1, Fondo per il Programma Nazionale di 
Ricerca e Progetti di Rilevante Interesse Nazionale (PRIN) funded by the European 
Commission under the NextGeneration EU programme, project ``Advanced optimization 
METhods for automated central veIn Sign detection in multiple sclerosis 
from magneTic resonAnce imaging (AMETISTA)'',  code: P2022J9SNP,
MUR D.D. financing decree n. 1379 of 1st September 2023 (CUP E53D23017980001),
project ``Numerical Optimization with Adaptive Accuracy and Applications 
to Machine Learning'',  code: 2022N3ZNAX MUR D.D. financing decree n. 973 of 30th
June 2023 (CUP B53D23012670006), and by Partenariato esteso FAIR ``Future Artificial 
Intelligence Research'' SPOKE 1 Human-Centered AI. Obiettivo 4, Project 
``Mathematical and Physical approaches to innovative Machine 
Learning technologies (MaPLe)'', Codice Identificativo EP\_FAIR\_002, CUP 
B93C23001750006.\\*[3mm]
The authors are indebted to an anonymous referee for a question that
led to a cleaner and more complete statement of the rate of
convergence results.

}%end footnotesize

\appendix
\renewcommand{\appendixname}{}
\section{Appendix. A variant of AS.5   }

We consider the following assumption in place of AS.5 and show the corresponding complexity result.

    {\bf Assumption AS.5bis:}
\textit{There exist the constants $\kappa_{Gg}\ge0$ and $\kappa_{Q}\ge0$ such that
\[
\Econd{k}{|(G_k - g_k)^Ts_k|}
\leq \kappa_{Gg}^2 \Econd{k}{\|s_k\|^2} + \kappa_Q Q_k,
\]
with $Q_k>0$, for all $k\geq 0$.
}

\llem{lemma:single-level-decrease_bis}{
Suppose that AS.3 and AS.5bis hold. Then, for $j\geq 0$,
\beqn{gen-decr-bis}
\Econd{j}{G_j^Ts_j}
\leq
-\frac{\tau\varsigma^2}{2\kappa_B}\,\Econd{j}{d_j^T\Delta_j}
+ \kappa_s^2 \left(\frac 1 2 \kappa_B+\kappa_{Gg}^2\right)\,\Econd{j}{\|\Delta_j\|^2}+ \kappa_Q Q_j,
\eeqn
where $\Delta_j\in\Re^n$ is the vector whose $i$-th component is $\Delta_{j,i}$.
}
\proof{Using AS.5bis and the second part of \req{inside-and-decr},
\[
\begin{array}{lcl}
\Econd{j}{G_j^T s_j}
& =& \Econd{j}{g_j^Ts_j}+\Econd{j}{(G_j-g_j)^Ts_j}\\*[1ex]
&\leq& \Econd{j}{g_j^Ts_j}+|\Econd{j}{(G_j-g_j)^Ts_j}|\\*[1ex]
&\leq& \Econd{j}{g_j^Ts_j}+\kappa_{Gg}^2|\Econd{j}{\|s_j\|^2}|+\kappa_Q Q_j\\*[1ex]
&\leq& \Econd{j}{g_j^Ts_j}+\kappa_{Gg}^2\kappa_s^2\Econd{j}{\|\Delta_j\|^2}+\kappa_Q Q_j\\
\end{array}
\]
and \req{gen-decr-bis} follows from \req{gen-decr0}.}

\lthm{theorem:convergence-bis}{
Suppose that AS.1--AS.4 and AS.5bis hold and that the
\algname\ algorithm is applied to problem \req{problem}. Then
\beqn{gradbound-bis}
\E{\average_{j\in\iiz{k}} \|d_j\|} \le \frac{\kap{conv}}{\sqrt{k+1}}+\kappa_Q
\frac{\sum_{j=0}^k Q_j}{\sqrt{k+1}},
\eeqn
where $\kap{conv}$ is given in 
\eqref{kapconv-def}.
}
\proof{The thesis follows proceeding as in the proof of Theorem \ref{theorem:convergence} replacing $\Gamma_0$ with $\Gamma_0+\kappa_Q \sum_{j=0}^k Q_j$.
}%epr

\noindent
If $\kappa_Q>0$ and $Q_j=\calO(1/j)$, we obtain that 
$\E{\average_{j\in\iiz{k}} \|d_j\|}$ decreases as  $\calO(\frac{\log(k)}{\sqrt {k+1}})$, while if $Q_j=\calO(1/j^p)$ with $p>\frac{1}{2}$, then  $\E{\average_{j\in\iiz{k}} \|d_j\|}$ decreases as  $\calO(\frac{1}{{k}^{p-1/2}})$.
Then, if $p=\frac{3}{4}$ the rate of convergence becomes the standard   $\calO(\frac{1}{{k}^4})$.

\end{document}